\documentclass[11pt,a4paper]{amsart}
\usepackage[utf8]{inputenc}
\usepackage[T1]{fontenc}
\usepackage[english]{babel}
\usepackage{amsmath}
\usepackage{amsfonts}
\usepackage{amssymb}
\usepackage{graphicx}
\usepackage{lmodern}
\usepackage[ left=3cm,right=3cm,top=3.5cm,bottom=5cm]{geometry} 
\usepackage{amsthm}
\usepackage{tikz}
\usepackage{verbatim}

\theoremstyle{plain}

\theoremstyle{definition}
\newtheorem{definition}{Definition}

\theoremstyle{remark}

\title[Calibration with CBO for interaction dynamics driven by NN]{Parameter calibration with Consensus-based Optimization  for interaction dynamics driven by neural networks} 

\author[S.~G\"ottlich]{Simone G\"ottlich}
\address[S.~G\"ottlich]{University of Mannheim}
\email[]{goettlich@uni-mannheim.de}

\author[C.~Totzeck]{Claudia Totzeck}
\address[C.~Totzeck]{University of Wuppertal}
\email[]{totzeck@uni-wuppertal.de}

\date{\today}

\begin{document}

\begin{abstract}We calibrate parameters of neural networks that model forces in interaction dynamics with the help of the Consensus-based global optimization method (CBO). We state the general framework of interaction particle systems driven by neural networks and test the proposed method with a real dataset from the ESIMAS traffic experiment. The resulting forces are compared to well-known physical interaction forces. Moreover, we compare the performance of the proposed calibration process to the one in \cite{NNinteraction} which uses a stochastic gradient descent algorithm.
\end{abstract}

\maketitle

\section{Introduction}
\label{sec:introduction}
Modelling interacting particle dynamics such as traffic, crowd dynamics, schools of fish and flocks of birds has attracted the attention of many research groups in the recent decades. Most models use physically-inspired interaction forces resulting from potentials to capture the observed behaviour. In fact, the gradient of the potential is used as driving force for interacting particle systems formulated with the help of ordinary differential equation (ODE). These models are able to represent the main features of the dynamics, but as for all models we cannot be sure that they deliver the whole truth. The idea in \cite{NNinteraction} was therefore to replace the physical-inspired models by neural networks, train the networks with real data and compare the resulting forces.

In the recent years it became obvious that neural networks are able to represent a lot of details from the dataset. It may be possible that there are details captured that are not even noticed by humans and therefore do not appear in physical models which are built to reproduce observations of the modeller.

In the following we recall the general dynamic of interaction particle systems driven by neural networks as proposed in \cite{NNinteraction}. Then we shortly describe the global optimization method 'Consensus-based optimization' that we use for the real-data based calibration the network. Finally, we present the numerical results obtained by the calibration process and compare them to the ones resulting from the calibration with the stochastic gradient descent method reported in \cite{NNinteraction}.

\section{Interacting particle systems driven by neural networks}
\label{sec:dynamics}

We consider interacting particle dynamics described by systems of ODEs of the following form
\begin{equation}\label{eq:ODE_totzeck}
	\frac{d}{dt} y_i =  \sum_{j=1}^N W^{i,j}_\theta(y_j -y_i), \quad y_i(0) = z_0^i,\quad i=1,\dots,N,
\end{equation}
where $W^{i,j}_\theta$ represents the interaction force resulting for $y_i$ in its interaction with $y_j.$ The initial condition of the particles is given by real dataset $z_0 = z(0).$ In order to compare the results to the ones in \cite{NNinteraction} we restrict the class of neural networks to feed-forward networks. However, note that the approach discussed here allows for general neural networks while the discussion in \cite{NNinteraction} considers feed-forward networks and can only be generalized to neural networks allowing for back propagation. 

\subsection{Feed-forward neural networks}
\label{subsec:NN_totzeck}
In the following we consider feed-forward artificial neural networks of the form
\begin{definition}\label{def:NN_totzeck}
	A \textit{feed-forward artificial neural network (NN)} is characterized by
	\begin{itemize}
		\item [-] Input layer: \vspace{-1.15em} $$a_1^{(1)} = 1, \quad a_k^{(1)} = x_{k-1},\; \text{ for }k\in \{ 2,\dots, n(1)+1\},$$  
		where $x \in \mathbb R^{n^{(1)}}$ is the input (feature) in \eqref{eq:ODE_totzeck} and $n^{(1)}$ is the number of neurons without the bias unit $a_1$.
		\item [-] Hidden layers: for $\ell \in \{ 2,\dots, L-1 \}, k \in \{2,\dots,n^{(\ell)} +1 \}$  \vspace{-0.85em} $$ a_1^{(\ell)} = 1,\quad a_k^{(\ell)} = g^{(\ell)}\left( \sum_{j=1}^{n^{(\ell -1)}+1} \theta_{j,k}^{(\ell-1)} a_j^{(\ell-1)} \right). $$  \vspace{-1.05em}
		\item [-] Output layer: \quad $a_k^{(L)} = g^{(L)} \left(\sum_{j=1}^{n^{(L -1)}+1} \theta_{j,k}^{(L-1)} a_j^{(L-1)}\right) \quad \text{for} \quad  k \in \{1,\dots, n^{(L)}\}$
	\end{itemize}
\end{definition}

Note that the output layer has no bias unit. The entry $\theta_{j,k}^{\ell}$ of the weight matrix $\theta^{(\ell)} \in \mathbb R^{n^{(\ell-1)} \times n^{(\ell)}}$ describes the weight from neuron $a_{j}^{(\ell-1)}$ to the neuron $a_k^{(\ell)}$. For notational convenience, we assemble all entries $\theta_{j,k}^{(\ell)}$ in a vector $\mathbb R^K$ with $$K := n^{(1)}\cdot n^{(2)} + n^{(2)} \cdot n^{(3)} + \dots + n^{(L-1)} \cdot n^{(L)}.$$ For the numerical experiment we use $g^{(\ell)}=\log(1+e^x)$ for $\ell =2,\dots,N-1$ and $g^{(L)}(x)=x.$
%An illustration of an NN with $L=3,$ one input and 5 units in the hidden layer can be found in Figure~\ref{fig:NN}. This is one of the networks used for the traffic application in the numerics section. We consider one lane, i.e.,~the difference of the positions of the cars is a scalar value. This difference $y_j -y_i$ is the input of the network. The force resulting from this interaction for $y_j$ is given by the output $a_1^{(3)}.$
For an illustration of the NN structure we refer the interested reader to \cite{NNinteraction}. In the numerical section we consider an NN with $L=3,$ one input and 5 units in the hidden layer.

\section{Parameter Calibration}
\label{sec:calibration}
We formulate the task of the parameter calibration as an optimization problem. Let $u \in \mathbb R^d$ denote the vector of parameters to be calibrated. This could be the weights of the neural network $\theta$ and some other parameters, as for example the average length $L$ and the maximal speed $v_{\max}$ of the cars which we will consider in the application. As we want the network to recover the forces hidden in the real data dynamics, we define the cost function for the parameter calibration as
\begin{equation}
	J(y,u) = \frac{1}{2} \int_0^T  \| y(t) - z(t) \|^2 dt + \frac{\delta}{2} | u - u_\text{ref} |^2,  
\end{equation}
where $z$ denotes the trajectories of the cars obtained by the traffic experiment, and $u_\text{ref}$ are reference values for the parameters. The parameter $\delta$ allows to balance the two terms in the cost functional. In case no reference values of the parameters are available, we set $\delta = 0$ in the numerical section.

\subsection{Consensus-based optimization (CBO)}
\label{subsec:cbo}
We solve the parameter calibration problem with the help of a Consensus-based optimization method \cite{CBO}. In more details, we choose the variant introduced in \cite{CBO_CarrilloJin} which is tailored for high-dimensional problems involving the calibration of neural networks.
The CBO dynamics is itself a stochastic interacting particle system with $N_\text{CBO}$ agents given by stochastic differential equations (SDEs). The evolution of the agents is influenced by two terms. On the one hand, there is a deterministic term that aims to confine the positions of the agents at a weighted mean. On the other hand, there is a stochastic term that allows for exploration of the state space. The details are the following
\begin{equation}
	du_t^i = -\lambda(u_t^i - v_f) dt + \sigma \text{diag}(u_t^i - v_f) dB_t^i, \quad i=1,\dots, N_\text{CBO}
\end{equation}	
with drift and diffusion parameters $\lambda, \sigma >0$, independent $d$-dimensional Brownian motions $B_t^i$ and initial conditions $u_0^i$ drawn uniformly from the parameter set of interest. A main role plays the weighed mean $$v_f = \frac{1}{\sum_{i=1}^{N_\text{CBO}} e^{-J(u_i)}} \sum_{i=1}^{N_\text{CBO}} u_i\, e^{-\alpha J(u_i)}.$$  %given by
%\[
%v_f = \frac{1}{\sum_{i=1}^{N_\text{CBO}} e^{-J(u_i)}} \sum_{i=1}^{N_\text{CBO}} u_i\, e^{-\alpha J(u_i)}.
%\]
By its construction, agents with lower cost have more weight in the mean as the ones with higher cost. The parameter $\alpha$ allows to adjust this difference of the weights. For more information on the CBO method and its proof of convergence on the mean-field level we refer the interested reader to \cite{TrendsCBO} and the references therein.
As indicated by the notation above, the agents used in the CBO method are different realizations of parameter vectors that we consider for the calibration. For the numerical results NN4 we consider a neural network with $13$ weights, i.e., $\theta \in \mathbb{R}^{13}$. Moreover, we assume the maximal speed $v_\text{max}$ as additional parameter. Hence, for fixed $t$ we have for the $i$-th CBO agent $u_t^i \in \mathbb R^{14}.$

\section{Numerical results and conclusion}
\label{sec:results}
For the calibration of the parameters we consider real data from the project ESIMAS \cite{dataTraffic}. As we want to compare the results to the well-known follow-the-leader model for traffic flow (LWR) we recall its details
\begin{subequations}\label{eq:LWR}
	\begin{align}
		\frac{d}{dt} y_i(t) &= f\left(\frac{y_{i+1}(t) -  y_i(t)}{L} \right), \quad i=1,\dots, N-1, \\
		\frac{d}{dt} y_N(t) &= v_\text{max}.
	\end{align} 
\end{subequations}
Here $f(\cdot)$ is either $v_\text{max}\log(\cdot)$ or $v_\text{max} (1 - 1/\cdot).$ To be prepared for a reasonable comparison, we consider for the neural network dynamics
\begin{subequations}\label{eq:NNvelo}
	\begin{align}
		\frac{d}{dt} y_i(t) &= W^{i,i+1}_\theta(y_{i+1}(t) -  y_j(t)), \quad i=1,\dots, N-1, \\
		\frac{d}{dt} y_N(t) &= v_\text{max}
	\end{align} 
\end{subequations}
supplemented with initial data $y(0) = z_0$. This leads to $u = (v_\text{max}, \theta).$ To evaluate the models and compute the corresponding cost we solve all ODEs with an explicit Euler scheme. For details we refer to \cite{NNinteraction}. The number in the notation $NN2$, $NN4$ and $NN10$ corresponds to the number of nonbias neurons in the hidden layer.

\subsection{Data processing and numerical schemes}
The data collection of the ESIMAS project contains vehicle data from 5 cameras that were placed in a $1 km$ tunnel section on the German motorway A3 nearby Frankfurt/Main \cite{dataTraffic}. %For the parameter estimation we extract sequences from the data that contain three or more vehicles in one lane. For simplicity, we restrict our considerations to the middle lane data and neglect the data of the $y$-coordinate. Thus, we have one-dimensional traffic data for the parameter estimation. 
The data is processed in the exact same way as in \cite{NNinteraction}. Files with the processed data can be found online\footnote{https://github.com/ctotzeck/NN-interaction}.

%First, we interpolate the position data that is supplemented with time stamps to a reference time discretization. Having all the data aligned to the reference time discretization, we filter sequences of data where two or more vehicles are present in the camera frame. This yields a database with various sequences of different length and with different number of vehicles that we use for the parameter identification. 
%
%Note that after this data extraction, the vectors containing the positions of the cars for each sequence are ordered. In fact, we pass this ordered vector to the neural network and hold on to the assumption that the interaction of the vehicles depends on the distance to the vehicle in front. This is why we have to prescribe some velocity for the first car in the data set, even in the case of the neural network without physical parameters, see \eqref{eq:NNvelo}.

The SDE which represents the CBO scheme is solved with the scheme proposed in \cite{CBO_CarrilloJin}. In particular, we set $dt = 0.05, \sigma_0 = 1, \lambda = 1$ and the maximal number of time steps to $100.$ The mini-batch size of the CBO scheme is $50$ and we have $100$ CBO agents in total. In each time step we update one randomly chosen mini-batch. The initial values are chosen as follows $$v_\text{max} \sim U([20,40]),\quad L \sim U([0,10]) \text{ and }\theta \sim U([-0.5,0.5]^K).$$

\subsection{Resulting forces and comparison}
\begin{figure}
	\includegraphics[scale=0.3]{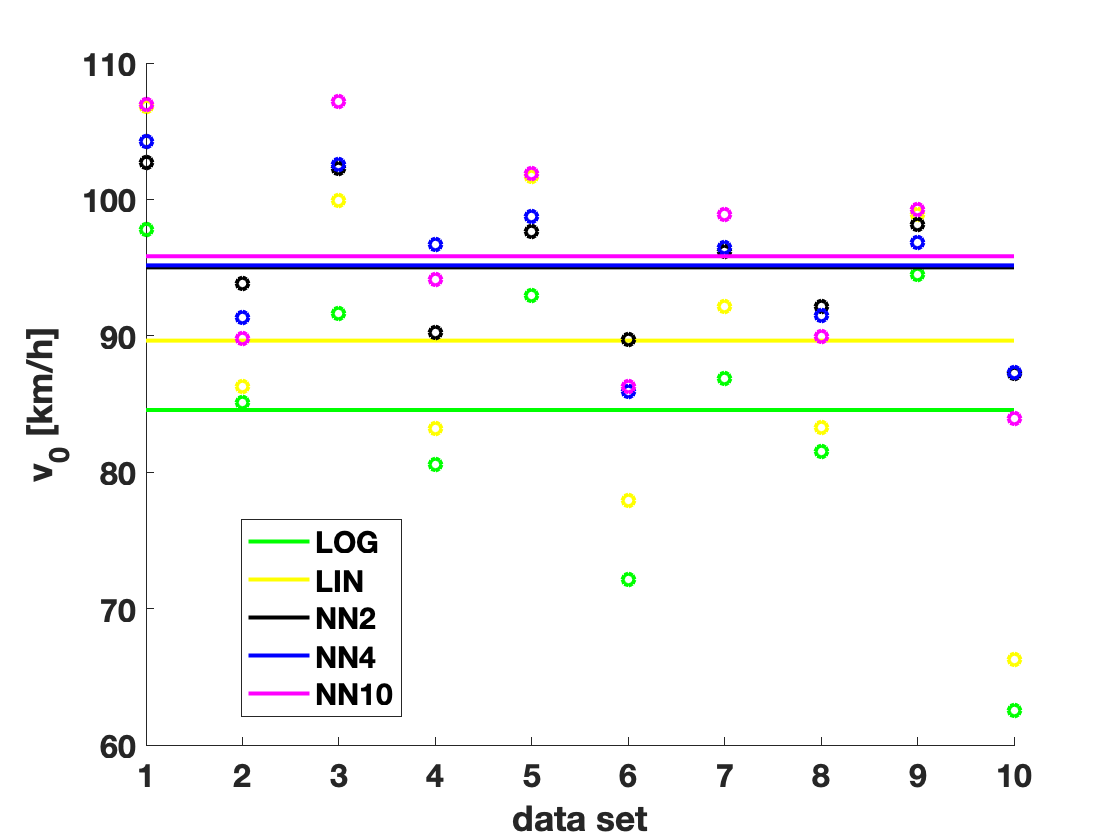}
	\includegraphics[scale=0.3]{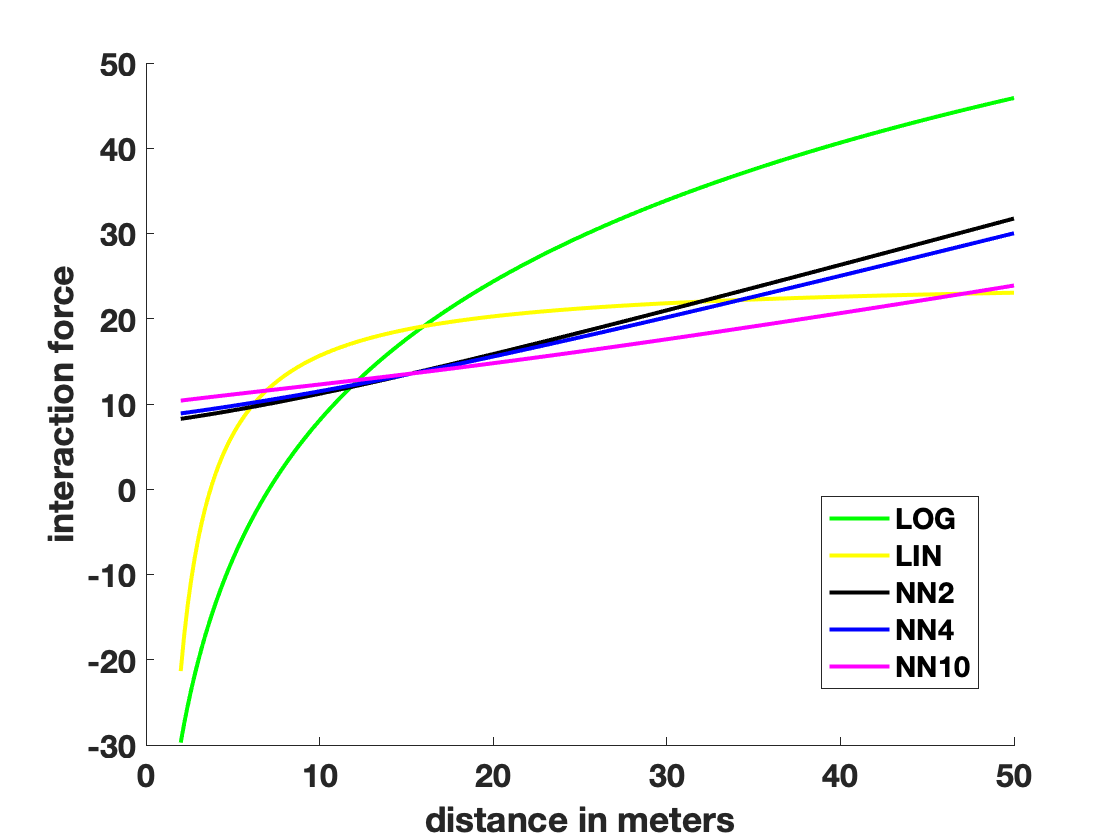}
	\caption{Average velocities and forces resulting from the parameter calibration and learning process. }
	\label{fig:velo_force}
\end{figure}

Figure~\ref{fig:velo_force} (left) shows the velocities resulting from the parameter calibration process. We find that the estimates velocities for the NN approaches are higher than the velocities of the LWR based models. The difference is most significant in data set $10.$ The plot on the right shows the average of the resulting forces for the different models. The forces of the NN approaches resemble linear approximations of the forces corresponding to the LWR models.
The car length $(L)$  appears only in the LWR models. Its optimized values for the different data sets are given in Table~\ref{tab:carlength}. We see that the lengths for the linear model are smaller than the ones in the logarithmic model. This is in agreement with the results obtained with stochastic gradient descent and shown in \cite{NNinteraction}. 
\begin{table}[ht!] \centering
	\begin{tabular}{l | c | c | c | c | c | c |c | c | c | c | c | c |}
		& 1 & 2 & 3 & 4 & 5 & 6 & 7 & 8 & 9 & 10 & average \\
		\hline
		Lin & 3.5969  &  3.76  &  4.17  &  2.19  &  3.02  &  2.81  &  5.92  &  5.86 &   2.14 &   3.65 & 3.71  \\
		\hline
		Log & 7.15 &   7.21  &  8.05  &  8.17  &  6.19  &  5.00  &  8.10  &  8.46  &  5.63   & 6.91 & 7.09  \\
		\hline 
	\end{tabular}
	\vspace{0.5em}
	\caption{ Car lengths (in $m$) estimated with the algorithm for the 10 data sets with the LWR-model with linear and logarithmic velocity.}
	\label{tab:carlength}
\end{table}
\begin{table}[ht!] \centering
	\begin{tabular}{l | c | c | c | c | c | c |c | c | c | c | c | c |}
		& 1 & 2 & 3 & 4 & 5 & 6 & 7 & 8 & 9 & 10 & average \\
		\hline
		NN2 & 47.95 &  46.49 &  98.07  & 44.97 &  23.69 &  29.72 &  40.69 &  55.75  & 11.50 &   68.91 & 46.77   \\
		\hline
		NN4 & 47.82 &  46.09 & 97.01 &  51.84 &  23.33 & 26.71 & 41.60 & \textbf{55.29} & 11.16 & 67.60 & 46.84 \\
		\hline
		NN10 & 47.90  & 45.78 &   99.20 & 42.50 & 22.16 & \textbf{24.40} & 41.18 &  56.68 &10.01 &  66.01 & 45.58 \\
		\hline
		Lin & \textbf{44.41} & \textbf{41.29}  & \textbf{93.73} &  \textbf{30.86} & \textbf{19.00} &  37.98  & \textbf{38.00} &   56.40 &  \textbf{8.18} &  \textbf{46.24} & \textbf{41.61}   \\
		\hline
		Log & 53.53 & 50.31 &109.36 & 65.24  & 26.50 & 52.93 & 38.09 & 58.22  & 14.54 &  52.75 & 52.15  \\
		\hline 
	\end{tabular}
	\vspace{0.5em}
	\caption{ Values of the cost functional estimated with the algorithm for the 10 data sets with the LWR-model with linear and logarithmic velocity and the three different neural network approaches}
	\label{tab:cost}
\end{table}
Finally, we summarize the cost values after parameter calibration in Table~\ref{tab:cost}. The least values of every column are highlighted. It is obvious that the LWR model with linear force outperforms the other models. The results of the NN approaches are better than the ones of the LWR model with logarithmic force.

\subsubsection{Comparison to calibration with stochastic gradient descent}
In comparison to the parameter calibration based on the stochastic gradient descent method reported in \cite{NNinteraction}, we find that the CBO approach finds better parameters for both LWR models. In fact, the resulting cost values are significantly smaller after the calibration with CBO. For the NN approaches the results are in good agreement. A clear decision in favour of the LWR approach or the NN ansatz was not possible based on the results of \cite{NNinteraction}. After the training with CBO the LWR with linear force seems to outperform  all other approaches. Note that we used NN with very simple structure here, it may be worth to test more sophisticated network structures in future work.

\bibliographystyle{plain}
\bibliography{biblio}
\end{document}